\documentclass[12pt,twoside]{article}
\usepackage{amsmath}
\usepackage{amsfonts}
\usepackage{amsthm}
\usepackage{amssymb}
\usepackage{graphicx}
\usepackage{multicol}

\begin{document}
\renewcommand{\baselinestretch}{1.5}
\title{{\normalsize
{\bf Ribbonness on classical link}}}
\author{{\footnotesize Akio Kawauchi}\\ 
\date{}
{\footnotesize{\it Osaka Central Advanced Mathematical Institute, Osaka Metropolitan University }}\\
{\footnotesize{\it Sugimoto, Sumiyoshi-ku, Osaka 558-8585, Japan}}\\
{\footnotesize{\it kawauchi@omu.ac.jp}}}
\maketitle
\vspace{0.25in}
\baselineskip=15pt
\thispagestyle{empty}
\newtheorem{Theorem}{Theorem}[section]
\newtheorem{Conjecture}[Theorem]{Conjecture}
\newtheorem{Lemma}[Theorem]{Lemma}
\newtheorem{Sublemma}[Theorem]{Sublemma}
\newtheorem{Proposition}[Theorem]{Proposition}
\newtheorem{Corollary}[Theorem]{Corollary}
\newtheorem{Claim}[Theorem]{Claim}
\newtheorem{Definition}[Theorem]{Definition}
\newtheorem{Example}[Theorem]{Example}

\begin{abstract} 
It is shown that if a link in 3-space bounds a proper oriented surface (without closed component) in the upper half 4-space, then the link bounds a proper oriented ribbon surface in the upper half 4-space 
which is a renewal embedding of the original surface. 
In particular, every slice knot is a ribbon knot, answering an old question by R. H. Fox affirmatively. 

\end{abstract}

\phantom{x}



\phantom{x}

\noindent{\bf 1. Introduction}

\phantom{x}

For a long time, the author has considered the (2,1)-cable of the figure-eight knot, which is not ribbon but rationally slice, as a candidate for a non-ribbon knot which might be slice (see \cite{Ka1, Ka2}). However, in 
\cite{DKMPS}, I. Dai, S. Kang, A. Mallick, J. Park and M. Stoffregen showed that it is not a slice knot. 
In this paper, the author comes back to elementary research beginning point of \cite{KSS} on the difference between a slice knot and a ribbon knot. 
Then it is concluded that every slice knot is a ribbon knot. 
More generally, it is shown that if a link in 3-space bounds a proper oriented surface (without closed component) in the upper half 4-space, then the link bounds a proper oriented ribbon surface in the upper half 4-space 
which is a renewal embedding of the original surface. 

This detailed explanation is done as follows. For a set $A$ in the 3-space 
${\mathbf R}^3=\{ (x,y,z)|\, -\infty<x,y,z<+\infty\}$ 
and an interval $J\subset {\mathbf R}$, let 
\[AJ=\{ (x,y,z,t)|\, (x,y,z)\in A,\, t\in J \}.\] 
The {\it upper-half 4-space} ${\mathbf R}^4_+$ is denoted by 
${\mathbf R}^3[0,+\infty)$. 
Let $k$ be a link in the 3-space ${\mathbf R}^3$, and $F$ a proper oriented surface in the 
upper-half 4-space ${\mathbf R}^4_+$ with $\partial F=k$. 
Let $b_j\, (j=1,2,\dots,m)$ be finitely many disjoint oriented bands spanning the link $k$ in ${\mathbf R}^3$, which are regarded as framed arcs spanning $k$ in ${\mathbf R}^3$.
Let $k'$ be a link in ${\mathbf R}^3$ obtained from $k$ by surgery along these bands. 
Then this band surgery operation is denoted by $k\to k'$. Let $k$ have $r$ knot components. 
If the link $k'$ has $r-m$ components, then the band surgery operation $k\to k'$ is called a {\it fusion}. 
If the link $k'$ has $r+m$ components, then the band surgery operation $k\to k'$ is called a {\it fission}. 
These terminologies are used in \cite{KSS}. 

A {\it band sum} $k\#_b o$ of a link $k$ and a trivial link $o$ of components $o_i\, (i=1,2,\dots,r)$ is a special fusion of the split sum $k+o$ along a disjoint band system $b_i\, (i=1,2,\dots,r)$ spanning $k$ and $o_i$ for every $i$. 
For the knot components $k_i\, (i=1,2,\dots,n)$ of $k$, assume that the band surgery operation $k\to k'$ induces the band surgery operation $k_i\to k'_i$ for all $i$. Then if the link $k'_i$ is a knot for all $i$, then 
the band surgery operation $k\to k'$ is called a {\it genus addition}. 

Every band surgery operation $k\to k'$ along a band system $b$ is realized as a proper surface $F_s^u$ in ${\mathbf R}^3[s,u]$ for any interval $[s,u]$, 
as follows (see \cite{KSS}):
\[
F_s^u\cap {\mathbf R}^3[t]=\left\{
\begin{array}{rl} 
k'[t],& \quad \mbox{for $\frac{s+u}{2} <t\leq u$},\\
(k\cup b)[t],& \quad \mbox{for $t=\frac{s+u}{2}$},\\
k[t],& \quad \mbox{for $s \leq t < \frac{s+u}{2}$}. 
\end{array}\right.
\] 
For every band surgery sequence $k_1\to k_2\to \dots \to k_{n-1} \to k_n$, 
the {\it realizing surface} $F_s^u$ in ${\mathbf R}^3[s,t]$ is given by the 
union 
\[F_{s_0}^{s_1}\cup F_{s_1}^{s_2}\cup \dots \cup F_{s_{m-2}}^{s_{m-1}}\cup F_{s_{m-1}}^{s_m}\]
for any division 
\[s=s_0<s_1<s_2<\dots<s_{m-1}<s_m=u\]
of the interval $[s,u]$. Note that the realizing surface 
$F_s^u$ in ${\mathbf R}^3[s,t]$ is uniquely determined up to smooth isotopies of ${\mathbf R}^3[s,t]$ keeping ${\mathbf R}^3[s]\cup {\mathbf R}^3[t]$ fixed. 
For a band surgery sequence $k_1\to k_2\to \dots \to k_{n-1} \to k_n$ 
where $k_1$ is a split sum $k'_1+o$ for a link $k'_1$ and a trivial link $o$ and $k_n$ is a trivial link $o'$, a {\it semi-closed realizing surface} 
$\mbox{scl}(F_s^u)$ in ${\mathbf R}^3[s,t]$ 
bounded by the link $k'_1$ in ${\mathbf R}^3$ is constructed as follows. 
\[\mbox{scl}(F_s^u)=F_s^u\cup d[s]\cup d'[u]\]
for disk systems $d, d'$ in ${\mathbf R}^3$ with $\partial d=o$ and 
$\partial d'=o'$. 
A {\it modified semi-closed realizing surface} $\mbox{scl}(F_s^u)^+$ 
of the band surgery sequence 
$k_1=k'_1+o\to k_2\to \dots \to k_{n-1} \to k_n=o'$ 
is a proper surface in ${\mathbf R}^3[s,+\infty)$ bounded by the link 
$k'_1$ obtained from $\mbox{scl}(F_s^u)$ by raising the level $s$ of the disk $d$ into the level 
$d+\varepsilon$ for a sufficiently small $\varepsilon>0$. 

Let $F$ be an $r$-component proper surface without closed 
component in the upper-half 4-space ${\mathbf R}^4_+$ which bounds 
a link $k$ in ${\mathbf R}^3$. 
By \cite{KSS}, the proper surface $F$ in ${\mathbf R}^4_+$ is equivalent to a modified semi-closed realizing surface
$\mbox{scl}(F^1_0)^+$ of a band surgery $k +o\to o'$ in ${\mathbf R}^4_+$.
Since the band system used for $k +o\to o'$ is made disjoint, 
the modified semi-closed realizing surface
$\mbox{scl}(F^1_0)^+$ is further equivalent to a modified semi-closed realizing surface $\mbox{scl}(F^1_0)^+$ 
of a band surgery sequence
\[ \mbox{(*)}\qquad\qquad k +o \to k_1\cup o \to k_2 \cup o \to k_3 \to o_4=o',\]
where 

\medskip 

\noindent{(0)} $k_1$ is a link of $r$ components and the operation 
$k +o \to k_1\cup o$ is a fusion fixing $o$, 

\medskip 

\noindent{(1)} the operation 
$k_1\cup o \to k_2 \cup o$ is a genus addition fixing $o$,

\medskip 

\noindent{(2)} the operation $k_2 \cup o \to k_3$ is a fusion along 
a band system connecting every component of $o$ to $k_2$ so that 
$k_3$ is a link with $r$ components,

\medskip 

\noindent{(3)} the operation $k_3 \to o_4=o'$ is a fission (cf. \cite{KSS}).

\phantom{x}

In particular, in the band surgery sequence (*) above, if the trivial link $o$ is taken the empty set $\emptyset$, then 
the step (2) is omitted and we have $k_2=k_3$. 
A proper surface $F$ in ${\mathbf R}^4_+$ is said to be 
{\it ribbon} if it is equivalent to a semi-closed realizing surface of a band surgery sequence (*) with $o=\emptyset$.

The purpose of this paper is to show the following theorem. 

\phantom{x}

\noindent{\bf Theorem~1.1.} Assume that a link $k$ in the 3-space 
${\mathbf R}^3$ bounds 
a proper oriented surface $F$ without closed component 
in the upper-half 4-space ${\mathbf R}^4_+$. Then 
the link $k$ in ${\mathbf R}^3$ bounds a ribbon surface $F'$ in 
${\mathbf R}^4_+$ which is a renewal embedding of $F$.

\phantom{x}

For a link $k$ in ${\mathbf R}^3$, let $g^*(k)$ be the minimal genus of a smoothly 
embedded connected proper surface in ${\mathbf R}^4_+$ bounded by $k$, and $g^*_r(k)$ the minimal genus of a connected ribbon surface in 
${\mathbf R}^4_+$ bounded by $k$. The following corollary is a direct consequence of Theorem~1.1. 

\phantom{x}

\noindent{\bf Corollary~1.2.} $g^*(k)=g^*_r(k)$ for every link $k$.

\phantom{x}

Since a slice knot in ${\mathbf R}^3$ is the boundary knot of a smoothly embedded proper disk in ${\mathbf R}^4_+$ and a ribbon knot in ${\mathbf R}^3$ is the boundary knot of a ribbon disk in ${\mathbf R}^4_+$,
Corollary~1.2 contains an affirmative answer to Fox Problem~25 \cite{Fox}. 

\phantom{x}

\noindent{\bf Corollary~1.3.} Every slice knot is a ribbon knot.

\phantom{x}

\noindent{\bf 2. Proof of Theorem~1.1}

\phantom{x}

The following lemma is a starting point of the proof of Theorem~1.1.

\phantom{x}

\noindent{\bf Lemma~2.1.} 
For a knot $k$ in ${\mathbf R}^3$, 
assume that a band sum $o'=k\#_b o$ of $k$ and a trivial link $o$ is a trivial knot in 
${\mathbf R}^3$. 
Then the knot $k$ is a ribbon knot in ${\mathbf R}^3$. 

\phantom{x}

\noindent{\bf Proof of Lemma~2.1.} Let $-k^*$ be the reflected inverse knot of a knot $k$ 
in ${\mathbf R}^3$. Then the connected sum $(-k^*)\# k$ is a ribbon knot in ${\mathbf R}^3$ (see \cite{FM}).
Since the band sum $o'=k\#_b o$ is a trivial knot, 
the connected sum $(-k^*)\# (k\#_b o)$ obtained by locally tying $-k^*$ to a string of $k$ in $k\#_b o$
is equivalent to the knot $(-k^*)\# o'=-k^*$. On the other hand, the knot $(-k^*)\# (k\#_b o)$ is a ribbon knot 
because it is a band sum of the ribbon knot $(-k^*)\#k$ and the trivial link $o$. 
Thus, the knot $-k^*$ is a ribbon knot. Since the reflected inverse knot of a ribbon knot is a ribbon knot, the knot $k$ is a ribbon knot. 
This completes the proof of Lemma~2.1. 
$\square$

\phantom{x}

\noindent{\bf Remark~2.2.} A ribbon presentation of the connected sum $(-k^*)\# k$ for a knot $k$ in ${\mathbf R}^3$ can be obtained from 
the chord diagram of any given diagram $D(k)$ of $k$ by \cite{Ka3,Ka4,Ka5,Ka6}. In fact, by \cite{Ka6}, let $D$ be an inbound diagram of $D(k)$ 
(namely, an arc diagram obtained from $D(k)$ by removing 
an open arc not containing a crossing point) with the end points in the infinite region of the plane ${\mathbf R}^3$, and 
$C$ a chord diagram of $D$. The diagram obtained from the based loop system of $C$ by surgery 
along a band system thickening the chord system is a ribbon presentation of the connected sum $(-k^*)\# k$. This is because the connected sum 
$(-k^*)\# k$ is the middle cross-section of the spun knot $S(k)$ of $k$ in ${\mathbf R}^4$ and the chord diagram $C$ canonically represents 
the spun knot $S(k)$ as a ribbon $S^2$-knot (see \cite{Ka3,Ka6,Ya}). 

\phantom{x}

Lemma~2.1 is generalized as follows.

\phantom{x}

\noindent{\bf Lemma~2.3.} 
For a link $k$ of $n$ knot components in ${\mathbf R}^3$, 
assume that a band sum $k\#_b o$ of $k$ and a trivial link $o$ is a ribbon link in ${\mathbf R}^3$. 
Then the link $k$ is a ribbon link in ${\mathbf R}^3$. 

\phantom{x}

\noindent{\bf Proof of Lemma~2.3.} 
For the components $k_i\, (i=1,2,\dots,n)$ of $k$, the band sum $k'=k\#_b o$ is the union of band sums $k'_i=k_i\#_b o_i\, (i=1,2,\dots,n)$. 
Let $o_{ij}\, (j=1,2,\dots,n_i)$ be the components of the trivial link $o_i$, and $b_{ij}$ the band spanning $k_i$ and $o_{ij}$
used for the band sum $k'_i=k_i\#_b o_i$ for all $j\, (j=1,2,\dots,n_i)$. 
Since the link $k'$ is a ribbon link with components $k'_i\, (i=1,2,\dots,n)$, there is a fusion $o'\to k'$ with a trivial link $o'$ consisting of 
fusions $o'_i\to k'_i\, (i=1,2,\dots,n)$. 
Let $o'_{ih}\, (h=1,2,\dots,m_i)$ be the components of $o'_i$, and $b'_{ih}\, (h=1,2,\dots,m_i)$ the bands used for the fusion $o'_i\to k'_i$. 
By band slides and by regarding bands as framed arcs, the bands 
$b_{ij} (i=1,2,\dots, n;\, j=1,2,\dots,n_i)$, $b'_{ih}\, (i=1,2,\dots, n;\, h=1,2,\dots,m_i)$ are made disjoint. 
Further, the bands $b_{ij}\, (j=1,2,\dots,n_i)$ are taken to be attached only to the component $o'_{i1}$. 
Let $B'_{ij}\, (i=1,2,\dots,n; j=1,2,\dots, m_i)$ be disjoint 3-balls in ${\mathbf R}^3$ containing the component $o'_{ij}$ in the interior. 
Let $d_{ij}\, (i=1,2,\dots,n;\, j=1,2,\dots,n_i)$ be a disjoint disk system bounded by the trivial loop system $o_{ij}\, (i=1,2,\dots,n;\, j=1,2,\dots,n_i)$ in ${\mathbf R}^3$. Let
$a'_{ih}\, (i=1,2,\dots, n;\, h=1,2,\dots,m_i)$ be a core arc system of the band system $b'_{ih}\, (i=1,2,\dots, n;\, h=1,2,\dots,m_i)$, 
and $a''_{ih}\, (i=1,2,\dots, n;\, h=1,2,\dots,m_i)$ an arc system obtained 
from $a'_{ih}\, (i=1,2,\dots, n;\, h=1,2,\dots,m_i)$ by deforming not to meet the disjoint disk system $d_{ij}\, (i=1,2,\dots,n;\, j=1,2,\dots,n_i)$. 
The deformation should be taken so that the arc system $a''_{ih}\, (i=1,2,\dots, n;\, h=1,2,\dots,m_i)$ is isotopic to the arc system $a'_{ih}\, (i=1,2,\dots, n;\, h=1,2,\dots,m_i)$ when the disk system 
$d_{ij}\, (i=1,2,\dots,n;\, j=1,2,\dots,n_i)$ is forgotten. 
Let $b''_{ih}\, (i=1,2,\dots, n;\, h=1,2,\dots,m_i)$ be the band system thickening the core arc system $a''_{ih}\, (i=1,2,\dots, n;\, h=1,2,\dots,m_i)$. 
Then the disjoint disk system $d_{ij}\, (i=1,2,\dots,n;\, j=1,2,\dots,n_i)$ can be moved into $B'_{i1}$ while keeping the band system $b''_{ih}\, (i=1,2,\dots, n;\, h=1,2,\dots,m_i)$ fixed. In this move, some parts of the band system 
$b_{ij}\, (i=1,2,\dots, n;\, j=1,2,\dots,n_i)$ may be moved.
Since $o_{i1}$ and $d_{ij}\, ( j=1,2,\dots,n_i)$ are disjoint 
except for the meeting part of the band system $b_{ij}\, ( j=1,2,\dots,n_i)$, 
there is a knot $k''_i$ such that the trivial knot $o_{i1}$ is the band sum 
$k''_i\#_b o_i$ using the bands $b_{ij}\, (j=1,2,\dots,n_i)$. 
By Lemma~2.1, the knot $k''_i$ is a ribbon knot and thus there is a fusion $o''_i\to k''_i$ for a trivial link $o''_i$ in ${\mathbf R}^3$. 
Note that the knot $k''_i$ is disjoint from $B'_{ij}\, (i=2,3,\dots,m_i)$, so that the trivial link $o''_i$ is movable into $B'_{i1}$ although some parts of the bands used for the fusion $o''\to k''$ may not be in $B'_{i1}$. 
The link $k$ is a fusion of the trivial link consisting of the split sum of $o'_i, \,(i=,2,3,\dots,n);\, o''_i\, (i=1,2,\dots,n)$, meaning that the link $k$ is a ribbon link.
This completes the proof of Lemma~2.3. 
$\square$

\phantom{x}

The proof of Theorem~1.1 is done as follows. 

\phantom{x}

\noindent{\bf Proof of Theorem~1.1.} 
Consider that a proper oriented surface $F$ is given by the sequence 
\[ k +o \to k_1\cup o \to k_2 \cup o \to k_3,\]
which are given by the band surgery operations that $k_3 \to k_2\cup o$ is a fission, $k_2\cup o\to k_1\cup o$ is a genus addition fixing $o$ and $k_1\cup o\to k+o$ is a fission fixing $o$, forming the inverse sequence 
\[k_3 \to k_2\cup o\to k_1\cup o\to k+o\]
of the sequence $k +o \to k_1\cup o \to k_2 \cup o \to k_3$. 
By band slides, note that  the link 
$k_2\cup o$  can be the split link $k_2+ o$. 
Replace the bands used for the genus addition 
$k_2\cup o\to k_1\cup o$ and the fission   $k_1 \cup o \to k+o$ 
by bands such that 

\medskip 

\noindent{(i)} every band does not change the attaching parts, and 

\medskip 

\noindent{(ii)} every band does not pass the trivial link $o$, and

\medskip 

\noindent{(iii)} every band is deformable into the original band if the trivial link $o$ is forgotten. 

\phantom{x}

Then the genus addition 
$k_2\cup o\to k_1\cup o$ changes into a genus addition 
$k_2+o\to k_1+o$ fixing $o$  and the fission 
$k_1 \cup o \to k+o$ changes into a fission 
$k_1 + o \to k+o$ fixing $o$, respectively, so that  
the sequence 
\[ k +o \to k_1\cup o \to k_2 \cup o \to k_3\]
changes into a sequence 
\[ k +o \to k_1+ o \to k_2 + o \to k_3,\]
where 
the operation $k+o \to k_1+ o$ is a fusion fixing $o$, the operation 
$k_1+o \to k_2 + o$ is a genus addition fixing $o$, 
and the operation $k_2 + o \to k_3$ is a fusion meaning that $k_3$ is a band sum $k_2\#_b o$ of $k_2$ and $o$. Since $k_3$ is a ribbon link, $k_2$ is a ribbon link by Lemma~2.3. Thus, there is a sequence  
\[ k \to k_1 \to k_2 \to o'_3,\]
where the operation $k \to k_1$ is a fusion, the operation 
$k_1 \to k_2 $ is a genus addition 
and the operation $k_2 \to o'_3$ is a fission with $o'_3$ a trivial link. 
This means that the link $k$ in ${\mathbf R}^3$ bounds a ribbon surface $F'$ in 
${\mathbf R}^4_+$ which is a renewal embedding of $F$.
This completes the proof of Theorem~1.1. 
$\square$

\phantom{x}

\noindent{\bf Acknowledgements.} 
This paper was completed during the author's stay at Dalian University of Technology, China from July 2, 2023 to July 21, 2023. 
The author would like to thank Feng Chun Lei and Fengling Li 
for their kind hospitalities. 
This work was partly supported  by JSPS KAKENHI Grant Numbers JP19H01788, JP21H00978  and MEXT Promotion of Distinctive Joint Research Center Program JPMXP0723833165.

\phantom{x}

\end{document}